\documentclass[12pt]{article}
\usepackage[utf8]{inputenc}

\usepackage{type1cm}         

\usepackage{makeidx}         
\usepackage{graphicx}        
\usepackage{multicol}        
\usepackage[bottom]{footmisc}

\usepackage{newtxtext}       %
\usepackage{newtxmath}       

\usepackage[a4paper, total={6in, 8in}]{geometry}

\usepackage{tikz}
\usetikzlibrary{plotmarks}
\usetikzlibrary{decorations.pathreplacing}
\usepackage{pgfplots}
\pgfplotsset{width=10cm,compat=1.9}
\usepgfplotslibrary{external}
\usepackage{soul,xcolor} 
\setstcolor{blue}       
\newcommand{\tcb}[1]{{\textcolor{blue}{#1}}} 

\makeatletter
\renewcommand{\fnum@figure}[1]{\textbf{\figurename~\thefigure}. }
\renewcommand{\fnum@table}[1]{\textbf{\tablename~\thetable}. }
\makeatother

\allowdisplaybreaks

\usepackage{hyperref}

\newcommand{\eq}[1]{(\ref{eq:#1})}

\newcommand{\rem}[1]{Remark~\ref{rem:#1}}
\newcommand{\prop}[1]{Proposition~\ref{prop:#1}}
\newcommand{\cor}[1]{Corollary~\ref{cor:#1}}

\newtheorem{Theorem}{Theorem}

\newtheorem{Proposition}{Proposition}
\newtheorem{Corollary}{Corollary}
\newtheorem{Remark}{Remark}

\DeclareMathOperator{\E}{\mathbb{E}}
\DeclareMathOperator{\Pb}{\mathbb{P}}

\newcommand{\TAU}{\tau}

\title{On fluctuation-theoretic decompositions \\via Lindley-type recursions}

\author{Onno Boxma\thanks{EURANDOM and Department of Mathematics and Computer Science, Eindhoven University of Technology, the Netherlands. \texttt{o.j.boxma@tue.nl}}
\thanks{Partly funded by the NWO Gravitation project N{\sc etworks}, grant number 024.002.003.} ,
Offer Kella\thanks{Department of Statistics, The Hebrew University of Jerusalem; Jerusalem 9190501, Israel. \texttt{offer.kella@gmail.com}} 
\thanks{Partially supported the Vigevani Chair in Statistics.} \ and Michel Mandjes\thanks{Korteweg-de Vries Institute for Mathematics, University of Amsterdam, Science Park 904, 1098 XH Amsterdam, The Netherlands. \texttt{m.r.h.mandjes@uva.nl}} \footnotemark[2]
}

\date{\today}

\begin{document}

\maketitle

\begin{abstract}
\noindent
Consider a L\'evy process $Y(t)$ over an exponentially distributed time $T_\beta$ with mean $1/\beta$.
We study the joint distribution of the running maximum $\bar{Y}(T_\beta)$ and the time epoch $G(T_\beta$) at which this maximum last occurs.
Our main result is a fluctuation-theoretic distributional equality: the vector ($\bar{Y}(T_\beta),G(T_\beta)$) can be written as a sum of two independent vectors, the first one being ($\bar{Y}(T_{\beta+\omega}),G(T_{\beta+\omega})$)
and the second one being the running maximum and corresponding time epoch under the restriction that the L\'evy process is only observed at Poisson($\omega$) inspection epochs (until $T_\beta$).
We first provide an analytic proof for this remarkable decomposition, and then a more elementary proof that gives insight into the occurrence of the decomposition
and into the fact that $\omega$ only appears in the right hand side of the decomposition. The proof technique underlying the more elementary derivation also leads to further generalizations of the decomposition,
and to some fundamental insights into a generalization of the well known Lindley recursion.
\end{abstract}


\bigskip
\noindent {\bf AMS Subject Classification (MSC2020):} Primary 60G50, 60G51; Secondary 60K25.

\newpage
\section{Introduction}
\label{intro}
Consider a general real-valued L\'evy process $Y \equiv \{Y(t), \, t \geq 0\}$, inspected at independent Poisson($\omega$) epochs $I_1,I_2,\ldots$ 
In \cite{bm21} a relation is established between (i) the distribution of the running maximum $\bar{Y}(T_\beta)$ of $Y$ under continuous observation
until an, independently sampled, exp($\beta$) distributed `killing time' $T_\beta$
and  (ii) the running maximum $Y_{\beta,\omega}$ of $Y$ until $T_\beta$ when restricting oneself to those Poisson($\omega$) inspection epochs. 
This also covers the infinite horizon case by taking $\beta \downarrow 0$, when imposing the
additional assumption that $Y$ has a negative drift.
The main result in \cite{bm21} is the following fluctuation-theoretic decomposition (with, throughout the entire paper, the symbol `$\sim$' denoting equality in distribution).
\label{sec:intro}
\begin{Theorem}
\label{Thm1}
The following distributional equality holds:
\begin{equation}
\bar Y(T_\beta) \sim \bar Y(T_{\beta+\omega})+Y_{\beta,\omega},
\label{onedimdeco}
\end{equation}
with the two random variables in the right hand side being independent.
\end{Theorem}

The practical motivation to study such relations between processes under continuous observation and under Poisson inspection
is that real-life processes usually are not observed continuously.
Examples abound in, {\em e.g.}, healthcare, reliability and insurance. In the latter area, 
one has the distinction between the ruin probability in the classical Cram\'er-Lundberg insurance risk model,
and the bankruptcy probability for the same model, as introduced in \cite{AGS}: bankruptcy occurs if the insurance company has a negative capital {\em at a Poisson inspection epoch}.
More specifically,  the ruin probability $p(u,T_\beta)$ before $T_\beta$ of an insurance company with initial capital $u$ can be written as ${\mathbb P}(\bar{Y}(T_\beta)>u)$,
whereas the (smaller) bankruptcy probability is given by ${\mathbb P}(Y_{\beta,\omega}>u)$.
The latter probability was obtained in \cite{AL} for $\beta=0$ and $Y$ a compound Poisson process with exponentially distributed claim sizes;
the case of generally distributed claim sizes was treated in \cite{BEJ}.
In reliability theory, similar quantities are highly relevant; think of the probability that the condition of a device (car, {\sc mri} scanner) exceeds a particular critical value at which, {\em e.g.}, repair is needed.

\smallskip

In \cite{bm21}, Theorem~\ref{Thm1} was proven relying on Wiener-Hopf factorizations for L\'evy processes and various manipulations of integrals.
Theorem~\ref{Thm1} indeed has the flavor of a Wiener-Hopf decomposition, but is actually of a quite different nature.
A famous result regarding the Wiener-Hopf decomposition of a L\'evy process $Y$ 
(see, {\em e.g.,} \cite[Thm.\ 45.7, p.\ 341]{sato}, \cite[Thm.\ 5, p.\ 159]{bertoin}, \cite[Thm.\ 6.16, p.\ 158]{kyprianou}) states that the random variables 
$\bar{Y}(T_\beta)$ and $Y(T_\beta) - \bar{Y}(T_\beta)$, which evidently sum to $Y(T_\beta)$, are independent.
An intrinsic difference with the setting of Theorem~\ref{Thm1}, is that there the two independent components in the sum in the right hand side add up to the left hand side only `in a distributional sense'.

\smallskip

The proof of Theorem~\ref{Thm1} in \cite{bm21} neither provides any intuition on the remarkably simple decomposition \eqref{onedimdeco}.
In addition, it did not offer any explanation for the intriguing fact that $\omega$ only appears in its right hand side; while clearly $\bar Y(T_{\beta+\omega})$ decreases in $\omega$ and $Y_{\beta,\omega}$ increases in $\omega$, they apparently do so in a way that these effects cancel. 
These observations have motivated us to seek a deeper insight into the decomposition in Theorem~\ref{Thm1};
further theoretical motivation was provided to us by several beautiful identities that were presented in \cite{ai17}, relating exit problems for L\'evy processes under permanent observation and their counterparts under Poisson inspections.
Theorem~\ref{Thm1} raised the following concrete questions with us, which we aim to answer in the present paper.
\begin{enumerate}
\item
Can Theorem~\ref{Thm1} be extended to the two-dimensional case, with as second component  the {\em time epochs} at which the respective running maxima occur?
\item
Is there a (relatively) elementary way to understand and prove Theorem~\ref{Thm1}?
\item
Are there any generalizations of Theorem~\ref{Thm1} in the one-dimensional case?
\item
Can it be explained that $\omega$ only appears in the right hand side of (\ref{onedimdeco})?
\end{enumerate}
The main results of the present paper are the affirmative answers to these four questions, including an explicit probabilistic decomposition
of the vector
($\bar{Y}(T_\beta),G(T_\beta)$), with $G(T_\beta$) the time epoch of the last occurrence of the running maximum,
into two independent vectors (and a direct way to derive it). The first of these vectors is ($\bar{Y}(T_{\beta+\omega}),G(T_{\beta+\omega})$)
and the second one is the running maximum and corresponding time epoch under the restriction that the L\'evy process is only observed at Poisson($\omega$) inspection epochs, until $T_\beta$.
More specifically:
Let $H_m$ denote the increment of $Y$ between two consecutive inspection epochs $I_{m-1}$ and $I_m$ (with $I_0:=0$). Define $S_0:=0$ and, for $n \in \mathbb N$,
\[
S_n := \sum_{m=1}^n H_m,\]
and the corresponding running maximum process $\bar{S}_n := {\rm max}\{S_0,S_1,\dots,S_n\}$.
Obviously, the number $N_{\beta,\omega}$ of Poisson inspection epochs until the exponentially distributed $T_\beta$ is geometrically distributed, with success probability $q:=\beta/(\beta+\omega)$,
\begin{equation}
\Pb(N_{\beta,\omega}=n) = \left(1-q\right)^n q , ~~~n=0,1,\dots ~.
\label{eqgeom}
\end{equation}
Notice that $Y_{\beta,\omega} = \bar{S}_{N_{\beta,\omega}}$.
With $G_{N_{\beta,\omega}}$ the last time epoch (of an inspection) before $T_\beta$ that this maximum occurs, we shall prove the following two-dimensional decomposition that generalizes Theorem~\ref{Thm1}:
\begin{equation}
(\bar Y(T_\beta), G(T_\beta)) \sim (\bar Y(T_{\beta+\omega}), G(T_{\beta+\omega}))+(\bar S_{N_{\beta,\omega}}, G_{N_{\beta,\omega}}),
\label{2decom}
\end{equation}
with the two pairs in the right hand side being independent.
In addition, our study of the above decompositions leads us to several variants and a generalization of the well known Lindley recursion,
and subsequently to a series of results for such recursions (including the uniqueness of their solution). Since, in the literature on fluctuation theory for L\'evy processes, distributional properties of the pair $(\bar Y(T_\beta), G(T_\beta))$ have been analyzed in great detail, we can use \eqref{2decom} to obtain their counterparts for $(\bar S_{N_{\beta,\omega}}, G_{N_{\beta,\omega}})$. For instance, we can use the argumentation of \cite[Remark 3.4]{bm21} to find the covariance of $\bar S_{N_{\beta,\omega}}$ and $G_{N_{\beta,\omega}}$ in case the driving L\'evy process is spectrally one-sided.

\smallskip

The remainder of the paper is organized as follows.
In Section~\ref{sec:WH} we state and prove
the two-dimensional decomposition in (\ref{2decom}).
The proof technique uses the same -- rather heavy -- machinery as was used in \cite{bm21} to prove the above Theorem~\ref{Thm1}.
In Section~\ref{onedimrecur} we consider an elementary Lindley-type recursion that eventually, with the right choice of its components, yields a very simple proof
of Theorem~\ref{Thm1}. Section~\ref{Lindley} discusses the question whether a particular generalized Lindley-type equation has a unique solution.
In Section~\ref{twodimrecur} we outline how (\ref{2decom}) can be obtained using a fairly straightforward two-dimensional extension of the Lindley recursion of Section~\ref{onedimrecur}.
Section~\ref{implic} presents various explicit results for spectrally positive
and spectrally negative L\'evy processes (covering both transforms and moments),
and Section~\ref{sec:concl} contains concluding remarks and possible future research directions.

\smallskip

Some notation will be used throughout the paper. 
In what follows, {\em a.s.}, {\em iff}, $x\vee y$ and $x^+$ abbreviate {\em almost surely}, {\em if and only if}, $\max(x,y)$ and $\max(x,0)$, respectively. We use the symbol ${\mathbb N}$ to denote $\{0,1,2,\ldots\}.$ Also, $1_{\{A\}}$ is the indicator function of the event $A$.

\section{Decomposition: the two-dimensional case}
\label{sec:WH}
In this section we present a first proof of (\ref{2decom}), which is Theorem~\ref{Thm2} below.
This theorem is a two-dimensional decomposition that generalizes Theorem~\ref{Thm1} by including, as second component,
the last  {time epochs} at which the various maxima occur. We thus answer the first question listed in Section~\ref{intro}.

Before stating
and proving the theorem, we first mention a few key results regarding the maximum of a L\'evy process over an exponentially distributed interval, mainly taken from  \cite[Thm.\ 6.16]{kyprianou}.
$T_\zeta$ once more denotes an exp($\zeta$) distributed random variable (sampled independently from the L\'evy process under study), and
$G(T_\zeta)$ denotes the epoch at which the maximum $\bar{Y}(T_\zeta)$ is {\it last} attained.
Relying on Wiener-Hopf theory, it was proven that the Laplace-Stieltjes transform of the joint distribution of $\bar{Y}(T_\zeta)$ and $G(T_\zeta)$ is given by
\begin{equation}
\E {\rm e}^{-\alpha \bar{Y}(T_\zeta) - \gamma G(T_\zeta)} = \frac{\kappa(\zeta,0)}{\kappa(\gamma + \zeta,\alpha)},
\label{WH1}
\end{equation}
where
\begin{equation}
\kappa(a,b) :=  k ~ {\rm exp}\left[\int_0^{\infty} \int_{(0,\infty)} \frac{1}{t}\, ({\rm e}^{-t} - {\rm e}^{-at-bx})  \,\Pb(Y(t) \in {\rm d}x) ~{\rm d}t\right],
\label{WH2}
\end{equation}
with $k$ some strictly positive constant.
As a consequence,
\begin{equation}
\E {\rm e}^{-\alpha \bar{Y}(T_\zeta) - \gamma G(T_\zeta)} = {\rm exp}\left[-\int_0^{\infty} \int_{(0,\infty)} \frac{1}{t}\, {\rm e}^{-\beta t}(1 - {\rm e}^{-\gamma t - \alpha x}) \,\Pb(Y(t) \in {\rm d}x) ~ {\rm d}t \right].
\label{WH3}
\end{equation}
It is readily verified that this in particular implies that
\begin{align*}
\E {\rm e}^{-\alpha \bar{Y}(T_\beta) - \gamma G(T_\beta)}
 = \:&\E {\rm e}^{-\alpha \bar{Y}(T_{\beta + \omega}) - \gamma G(T_{\beta + \omega})}\:\lambda(\alpha,\beta,\gamma,\omega),
\label{WH4}
\end{align*}
where
\begin{equation}
\label{lambda}
\lambda(\alpha,\beta,\gamma,\omega):={\rm exp}\left[-\int_0^{\infty} \int_{(0,\infty)} \frac{1}{t} {\rm e}^{-\beta t}(1 - {\rm e}^{-\omega t})(1 - {\rm e}^{-\gamma t - \alpha x}) ~\Pb(Y(t) \in {\rm d}x) ~{\rm d}t \right].
\end{equation}
We are now ready to formulate and prove the two-dimensional decomposition.
\begin{Theorem}
\label{Thm2}
The following distributional equality holds:
\[(\bar Y(T_\beta), G(T_\beta)) \sim (\bar Y(T_{\beta+\omega}), G(T_{\beta+\omega}))+(\bar S_{N_{\beta,\omega}}, G_{N_{\beta,\omega}}),\]
with the two pairs in the right hand side being independent.
\end{Theorem}
{\it Proof.}
Based on the above, it suffices to show that the Laplace-Stieltjes transform of the joint distribution of $(\bar{S}_{N_{\beta,\omega}},G_{N_{\beta,\omega}})$ is given by the expression for $\lambda(\alpha,\beta,\gamma,\omega)$ in  \eqref{lambda}.
Let $E_n(\beta + \omega)$ denote an Erlang distributed random variable with $n$ phases that are exponentially distributed with rate $\beta + \omega$. Then,
with $U_{N_{\beta,\omega}}$ the time of the $N_{\beta,\omega}$-th inspection epoch, recalling that $q=\beta/(\beta+\omega)$,
\[{\mathbb E}\,{\rm e}^{-\alpha S_{N_{\beta,\omega}}-\gamma U_{N_{\beta,\omega}}} = 
\sum_{n=0}^\infty \left(1-q\right)^n q \,{\mathbb E}\,{\rm e}^{-\alpha Y(E_n(\beta+\omega))-\gamma E_n(\beta+\omega)}.\]
With $\varphi(\alpha)$ the Laplace exponent of the L\'evy process $Y$, so that $\E {\rm e}^{\alpha Y(t)} = {\rm e}^{-\varphi(\alpha)t}$, we thus find
\[{\mathbb E}\,{\rm e}^{-\alpha Y(E_n(\beta+\omega))-\gamma E_n(\beta+\omega)}= \left(\frac{\beta+\omega}{\beta+\omega-\varphi(\alpha)+\gamma}\right)^n,\]
and hence
\begin{equation}\label{factor}{\mathbb E}\,{\rm e}^{-\alpha S_{N_{\beta,\omega}}-\gamma U_{N_{\beta,\omega}}} =\frac{\beta}{\beta+\omega} \frac{\beta+\omega-\varphi(\alpha)+\gamma}{\beta-\varphi(\alpha)+\gamma}.\end{equation}
The idea is now to find an alternative expression for the two factors in the right hand side in \eqref{factor}. To this end, we use the
Frullani integral: 
\[\frac{\beta}{\beta+\omega} = \exp\left[-\int_0^\infty\frac{1}{t}{\rm e}^{-\beta t} (1-{\rm e}^{-\omega t})\,{\rm d}t \right], \]
which is easily seen to be true by replacing ${\rm e}^{-\beta t}(1-{\rm e}^{-\omega t})/t$ by $\int_\beta^{\beta+\omega} {\rm e}^{-yt} {\rm d}y$.
It thus trivially follows that also
\[\frac{\beta}{\beta+\omega} = \exp\left[-\int_0^\infty \int_{\mathbb R}\frac{1}{t}{\rm e}^{-\beta t}(1-{\rm e}^{-\omega t})\,{\mathbb P}(Y(t)\in{\rm d}x)\,{\rm d}t\right].\]
Then we apply the same techniques to $(\beta -\varphi(\alpha) + \gamma)/(\beta+\omega -\varphi(\alpha)+\gamma)$, i.e., the second factor in the right hand side of \eqref{factor}. This leads to 
\begin{align}{\mathbb E}\,&{\rm e}^{-\alpha S_{N_{\beta,\omega}}-\gamma U_{N_{\beta,\omega}}}
\label{upperlowereq}
\\&= \exp\left[-\int_0^\infty \int_{\mathbb R}\frac{1}{t}{\rm e}^{-\beta t}(1-{\rm e}^{-\omega t})(1-{\rm e}^{-\gamma t-\alpha x})\,{\mathbb P}(Y(t)\in{\rm d}x)\,{\rm d}t\right].
\nonumber
\end{align}
The argumentation now follows that of the proof of \cite[Thm.\ 1]{WH_enc}, dealing with the Wiener-Hopf decomposition for random walks. More precisely, with the same line of reasoning as in the proof of part (i) of \cite[Thm.\ 1]{WH_enc}, it is obtained that $(\bar S_{N_{\beta,\omega}}, G_{N_{\beta,\omega}})$, being a geometric sum of i.i.d.\ random variables, is infinitely divisible, and in addition it follows that the random pairs (being in addition infinitely divisible)
\begin{equation}
    \label{2vect}
(\bar S_{N_{\beta,\omega}}, G_{N_{\beta,\omega}})\:\:\:\mbox{and}\:\:\:(S_{N_{\beta,\omega}} - \bar S_{N_{\beta,\omega}}, U_{N_{\beta,\omega}} - G_{N_{\beta,\omega}})\end{equation}
are independent. Then, as in the proof of part (ii) of \cite[Thm.\ 1]{WH_enc}, 
observing that $\bar{S}_{N_{\beta,\omega}}$ is nonnegative and $S_{N_{\beta,\omega}}-\bar{S}_{N_{\beta,\omega}}$ is nonpositive,
we conclude that
\begin{align}{\mathbb E}\,&{\rm e}^{-\alpha \bar S_{N_{\beta,\omega}}- \gamma G_{N_{\beta,\omega}}}
\label{uppereq}
\\&
= \exp\left[-\int_0^\infty \int_{(0,\infty)}\frac{1}{t}{\rm e}^{-\beta t}(1-{\rm e}^{-\omega t})(1-{\rm e}^{-\gamma t-\alpha x})\,{\mathbb P}(Y(t)\in{\rm d}x)\,{\rm d}t\right];
\nonumber
\end{align}
see also Remark \ref{R1} below regarding the second vector in \eqref{2vect}. This confirms 
that the Laplace-Stieltjes transform of the joint distribution of ($\bar S_{N_{\beta,\omega}},G_{N_{\beta,\omega}})$ is given by the expression for $\lambda(\alpha,\beta,\gamma,\omega)$ in \eqref{lambda}.
\hfill$\Box$

\begin{Remark}\label{R0}{\em 
Combining Equation \eqref{WH1} with the decomposition in Theorem~\ref{Thm2},
we conclude that 
the following expression holds for the joint transform of the running maximum (until $T_\beta$) at inspection epochs, and the last epoch at which that running maximum occurred:
\begin{equation}
{\mathbb E} {\rm e}^{-\alpha \bar S_{N_{\beta,\omega}} - \gamma G_{N_{\beta,\omega}}} = \frac{\kappa(\beta,0)}{\kappa(\beta+\gamma,\alpha)} \frac{\kappa(\beta+\omega+\gamma,\alpha)}{\kappa(\beta+\omega,0)} ,
\label{doubleLT}
\end{equation}
with $\kappa(\cdot,\cdot)$ as defined in \eqref{WH2}. }
\end{Remark}

\begin{Remark}\label{R1}\rm
It is noted that this Wiener-Hopf based argumentation,  that led to (\ref{uppereq}) from (\ref{upperlowereq}), also yields
\begin{align}{\mathbb E}\,&{\rm e}^{-\alpha (S_{N_{\beta,\omega}}-\bar S_{N_{\beta,\omega}})- \gamma (U_{N_{\beta,\omega}}-G_{N_{\beta,\omega}})}\nonumber
\\&
= \exp\left[-\int_0^\infty \int_{(-\infty,0)}\frac{1}{t}{\rm e}^{-\beta t}(1-{\rm e}^{-\omega t})(1-{\rm e}^{- \gamma t-\alpha x})\,{\mathbb P}(Y(t)\in{\rm d}x)\,{\rm d}t\right],
\label{lowereq}
\end{align}
with, as mentioned, the two pairs in \eqref{2vect} being independent. In addition, $S_{N_{\beta,\omega}}-\bar S_{N_{\beta,\omega}}$ is distributed as the corresponding running {\it minimum}. The transform appearing in \eqref{lowereq} can be expressed in terms of the function $\hat\kappa(\cdot,\cdot)$, as defined in \cite[Theorem 6.16]{kyprianou}; cf.\ \eqref{doubleLT}.
\hfill$\Diamond$
\end{Remark}
\begin{Remark}\rm
A related two-dimensional decomposition, restricted to the infinite-horizon case (corresponding to $\beta=0$), is given in \cite[Prop.\ 2]{ai17}.
A careful reading of the proof of that proposition reveals that one could possibly also apply the elegant ideas there to the finite horizon case. \hfill$\Diamond$
\end{Remark}

In Section~\ref{twodimrecur} we provide an alternative proof of Theorem~\ref{Thm2}, which we feel is more direct and intuitive. As a preparation for that, we first focus on its one-dimensional counterpart, as stated in Theorem~\ref{Thm1}.

\section{A one-dimensional recursion and its implications}
\label{onedimrecur}
In this section we aim to provide answers to Questions 2 and 3 posed in Section~\ref{sec:intro}: we wish to develop a more intuitive understanding of Theorem~\ref{Thm1}, and we wish to explore whether it can be further generalized (while remaining in the one-dimensional setting).

\medskip

Let $u_{n+1},v_n\ge 0$ for $n\ge 0$, assume $z_0=0$ and consider the following recursion between some deterministic quantities:
\begin{equation}\label{eq:Z}
z_{n+1}=u_{n+1}+(z_n-v_n)^+\,.
\end{equation}
To interpret this recursion, it may help the reader to think of a determinstically evolving queue in which there is service in order of arrival. Then $z_{n+1}$ corresponds to the sojourn time of the ($n+1$)-st customer who has service time $u_{n+1}$ and where the interarrival time since the $n$-th customer is $v_n$.

Noting that $z_1=u_1$, then with $w_n=z_{n+1}-u_{n+1}$ (think of the waiting time in the above-mentioned queue) and $x_n=u_n-v_n$ we have that the above recursion can be rewritten in a simpler form: $w_0=0$ and,
for every  $n=1,2,\dots$,
\begin{equation}\label{eq:W}
w_n=(w_{n-1}+x_n)^+\,.
\end{equation}
This is the well known Lindley recursion \cite{Lindley} of which the solution is
\begin{equation}\label{eq:W1}
w_n=s_n-\min_{k\in\{0,\ldots,n\}}s_k=\max_{k\in\{0,\ldots,n\}}(s_n-s_k)\,,
\end{equation}
where $s_0=0$ and $s_n=\sum_{i=1}^nx_i$ for $n=1,2,\dots$.

Recalling \eq{Z} and $w_n=z_{n+1}-u_{n+1}$, we therefore immediately have
\begin{align}\label{eq:ZW}
z_{n+1}&=u_{n+1}+w_n, \nonumber\\
w_n&=(z_n-v_n)^+\,,
\end{align}
where it should be noted that $w_n$ is a function of $u_1,\ldots,u_n,v_1,\ldots,v_n$, whereas $z_{n+1}$ is a function of $u_1,\ldots,u_{n+1},v_1,\ldots,v_n$. 

Turning to random variables, this seemingly elementary observation leads to some interesting conclusions.

\begin{Proposition}[Basic decomposition]\label{prop:ind}
Assume that $Z_0=0$ a.s.\ and that for each $n\in\{0,1,\ldots,K\}$ (where $K\in {\mathbb N}$ is finite or infinite),
\begin{equation}\label{eq:Z0}
Z_{n+1}\sim U_{n+1}+(Z_n-V_n)^+,
\end{equation}
where $U_{n+1},V_n,Z_n$ are assumed independent. Then for each $n\in\{0,1,\ldots,K\}$ we have that
\begin{align}\label{eq:ZW0}
Z_{n+1}&\sim U_{n+1}+W_n,\nonumber\\
W_n&\sim (Z_n-V_n)^+\,,
\end{align}
where $W_n\sim S_n-\min_{k\in\{0,\ldots,n\}}S_k$, $(S_k)_{k\in{\mathbb N}}$ is a random walk with independent (not necessarily identically distributed) increments
distributed like $U_i-V_i$ with $U_i,V_i$ independent, and where the pairs of random variables on the right of both equations in \eq{ZW0} are independent.
\end{Proposition}
{\it Proof}.
Take a sequence of independent random variables $U'_1,U'_2,\ldots$,
$V'_0,V'_1,\ldots$ such that
$U'_{n+1}\sim U_{n+1}$ and $V'_n\sim V_n$ for each $n\in{\mathbb N}$. Set $Z'_0=0$ and 
\begin{equation}\label{eq:Zprime}
 Z'_{n+1}=U'_{n+1}+(Z'_n-V'_n)^+.   
\end{equation}
It follows by induction that $Z_n\sim Z'_n$ for each $n\in{\mathbb N}$, and thus necessarily
\eqref{eq:ZW0} holds,
where (for $n\in{\mathbb N}$)  $W_n\sim S_n-\min_{k\in\{0,\ldots,n\}}S_k$, with $S_0=0$, $S_k=\sum_{i=1}^k(U'_i-V'_i)$ and $W_n$ is independent of $U_{n+1}$. Observe that the independence of $Z_n,V_n$ does not have to be proven as it was assumed.
\hfill$\Box$

\begin{Remark}\label{rem:iid}\rm
A special case of \prop{ind} is
\begin{equation}
Z_{n+1}\sim U+(Z_n-V)^+,
\end{equation}
where $U,V,Z_n$ are assumed independent. This means that in this case the increments of the random walk $S_n$ in the proof of the basic decomposition of \prop{ind} are i.i.d. As a consequence, $W_n\sim\max_{k\in\{0,\ldots,n\}}S_k$, which is a monotone sequence that converges in distribution to $W_\infty\sim \sup_{k\in {\mathbb N}}S_k$ as $n\to\infty$. Therefore $Z_n$ also converges in distribution and we have that
\begin{align}
Z_\infty&\sim U+W_\infty, \nonumber\\
W_\infty&\sim (Z_\infty-V)^+.
\label{ZWOinf}
\end{align}
We note that $W_\infty$ and $Z_\infty$ are either both a.s.\ finite or both a.s.\ infinite. As is well known, they are a.s.\ finite iff $S_n\to-\infty$ a.s. This, {\em e.g.,} follows from~\cite[Thm.~8.2.5, p.\ 264]{c74}. A sufficient condition is that ${\mathbb E}U<{\mathbb E}V$, which also includes the case that ${\mathbb E}U<\infty$ and ${\mathbb E}V=\infty$. \hfill$\Diamond$
\end{Remark}

\begin{Remark}\label{rem:iid1}\rm
For the case of \rem{iid}, let $N$ be a nonnegative integer-valued random variable (not necessarily a.s.\ finite). Denote $p_i:={\mathbb P}(N=i)$ for $i\in\{0,1,\ldots, \infty\}$. Let $W_N$ denote a random variable having the distribution
\begin{equation}
{\mathbb P}(W_N\in \cdot)=\sum_{0\le i\le \infty}p_i\,{\mathbb P}(W_i\in \cdot) ,
\end{equation}
and similarly for $Z_N$ and $Z_{N+1}$. Then we clearly have that
\begin{align}
Z_{N+1}&\sim U+W_N, \nonumber\\
W_N&\sim (Z_N-V)^+,
\label{ZWON}
\end{align}
where $W_N \sim \max_{k \in\{0,\ldots,N\}} S_k$ and
where the pairs of random variables appearing on the right of each of the two equations are assumed independent. 

We emphasize that if we would not have assumed that $U_{n+1}\sim U$ and $V_n\sim V$ for $n\in {\mathbb N}$, it would no longer be necessarily true that $U_{N+1},W_N$ are independent nor that $V_N,Z_N$ are independent, where $U_{N+1}$ and $V_N$ are defined in the obvious way. \hfill$\Diamond$
\end{Remark}
So far we have derived the very straightforward sets of recursions (\ref{eq:ZW0}) (for finite $n$), (\ref{ZWOinf}) (in steady state) and (\ref{ZWON}) (for random index $N$).
We shall now exploit these recursions via suitable choices of the input variables $U_{n+1},V_n$, to obtain an elementary proof of Theorem~\ref{Thm1}
as well as generalizations of that theorem.
In particular, we first relate $U_{n+1},V_n$ to values of a L\'evy process at Poisson inspection epochs.
The next three corollaries successively use the recursions (\ref{eq:ZW0}) from Proposition~\ref{prop:ind}, (\ref{ZWOinf}) from Remark~\ref{rem:iid} and (\ref{ZWON}) from Remark~\ref{rem:iid1}.

\begin{Corollary}\label{cor:Levy}
Assume that $Y\equiv \{Y(t),t\ge 0\}$ is a L\'evy process and that $\TAU_1,\TAU_2,\ldots$ are independent exponentially distributed random variables (not necessarily identically distributed) which are independent of $Y$.
Finally, denote $\bar Y(t):=\sup_{s\in[0,t]}Y(s)$, $\Sigma_0:=0$ and $\Sigma_n:=\sum_{i=1}^n \TAU_i$ for $n\in{\mathbb N}$. Then, for every $n\in{\mathbb N}$,
\begin{align}\label{eq:levy1}
\bar Y(\Sigma_{n+1}) &\sim U_{n+1}+W_n, \nonumber\\
W_n&\sim (\bar Y(\Sigma_n)-V_n)^+,
\end{align}
where $U_{n+1}\sim \bar Y(\TAU_{n+1})$, $V_n\sim \bar Y(\TAU_{n+1})-Y(\TAU_{n+1})$, $W_n=S_n-\min_{k\in\{0,\ldots, n\}}S_k$, $S_n$ is a random walk with increments distributed like $U_i-V_i$ and the pairs of random variables appearing on the right of each equation in \eq{levy1} are independent.
\end{Corollary}

{\it Proof.}
First introduce the `shift operator':  $\sigma_s Y(t):=Y(s+t)-Y(s)$ for $s,t\ge 0$. Then, for any $t, \TAU\ge 0$, we have that
\begin{equation}
\bar Y(\TAU+t)=\bar Y(\TAU) \vee (Y(\TAU)+\overline{\sigma_{\TAU} Y}(t))=\bar Y(\TAU)+(\overline{\sigma{_\TAU} Y}(t)-(\bar Y(\TAU)-Y(\TAU)))^+\,.
\label{eq22}
\end{equation}
By the stationary independent increments and strong Markov property of $Y$
we have that the random variable $\overline{\sigma_{\TAU} Y}(t)$ and random pair ($\bar Y(\TAU),\bar Y(\TAU)-Y(\TAU)$) are independent with $\overline{\sigma_{\TAU} Y}(t) \sim \bar Y(t)$. Therefore, from (\ref{eq22}), for each $t\ge 0$,
\begin{equation}
\label{eq20}
\bar Y(\TAU+t)\sim U+(\bar Y(t)-V)^+,
\end{equation}
with $U\sim \bar Y(\TAU)$ and $V\sim \bar Y(\TAU)-Y(\TAU)$; and $Y(t)$ is independent of ($U,V$).
Now assume in addition that $\TAU$ is exponentially distributed, independent of $Y$.
Then
the aforementioned Wiener-Hopf factorization for L\'evy processes
implies that $U$ and $V$ are independent.
Finally replace $\TAU$ by $\TAU_{n+1}$ and $t$ by $\Sigma_n$. Then, when $Z_n=\bar Y(\Sigma_n)$ and $U_{n+1},V_n$ are as assumed, \eq{Z0} is satisfied,
so the result is a direct consequence of \prop{ind}.\hfill$\Box$

\bigskip

We emphasize that, in \cor{Levy}, $U_i\sim \bar Y(\TAU_i)$ whereas $V_i\sim \bar Y(\TAU_{i+1})-Y(\TAU_{i+1})$ and thus, unless $\TAU_i$ and $\TAU_{i+1}$ are identically distributed,
it does {\em not} necessarily hold that $U_i-V_i \sim Y(\TAU_i)$.
However when $\TAU_1,\TAU_2,\ldots$ are i.i.d., then we indeed have, by the independence of $\bar Y(\TAU)$ and $\bar Y(\TAU)-Y(\TAU)$, that $U_i-V_i\sim Y(\TAU)$ and hence $S_k\sim Y(\Sigma_k)$ for every $k\in {\mathbb N}$.
 
The last paragraph, \cor{Levy} and \rem{iid} immediately imply the following, which may also be concluded from (the proof of) \cite[Prop.~1]{ai17}.
\begin{Corollary}\label{cor:Levy1}
If, in addition to the assumptions of \cor{Levy}, $\TAU_1,\TAU_2,\ldots$ are also identically distributed like $\TAU$, then
\begin{align}\label{eq:levy2}
\bar Y(\infty) &\sim U+W_\infty , \nonumber\\
W_\infty&\sim (\bar Y(\infty)-V)^+ ,
\end{align}
where $U\sim \bar Y(\TAU)$, $V\sim \bar Y(\TAU)-Y(\TAU)$, $\bar Y(\infty) \sim\sup_{t\ge 0}Y(t)$, $W_\infty\sim\sup_{n\in{\mathbb N}}Y(\Sigma_n)$, and the pairs of random variables on the right of each of the equations in \eq{levy2} are independent.
\end{Corollary}

Similarly, \cor{Levy}, \cor{Levy1} and \rem{iid1} also imply the following claim.
\begin{Corollary}\label{cor:Levy2}
With the assumptions of \cor{Levy1}, if $N$ is a nonnegative integer valued random variable (not necessarily a.s.\ finite) which is also independent of $Y$, then
\begin{align}\label{eq:levy3}
\bar Y(\Sigma_{N+1})&\sim U+W_N, \nonumber\\
W_N&\sim (\bar Y(\Sigma_N)-V)^+\,,
\end{align}
where $U,V$ are as defined in \cor{Levy1}, $W_N\sim\max_{k\in\{0,\ldots,N\}}Y(\Sigma_k)$, and the pairs of random variables on the right of each of the equations in \eq{levy3} are independent.
\end{Corollary}

Corollary \ref{cor:Levy3} below, being an immediate consequence of Corollary~\ref{cor:Levy2}, directly implies Theorem~\ref{Thm1}.
In the sequel, $M \sim {\mathbb G}(p)$ means that $M$ is geometrically distributed, in the sense that
${\mathbb P}(M=n) = p(1-p)^ {n-1}$, $n=1,2,\dots$ for $p\in[0,1]$. As a consequence, we have that $N_{\beta,\omega} +1 \sim {\mathbb G}(q)$, recalling that $q=\beta/(\beta+\omega)$, cf.\ (\ref{eqgeom}). In addition, ${\mathbb B}(q)$ corresponds to the Bernoulli distribution with success probability $q\in[0,1]$.
\begin{Corollary}\label{cor:Levy3}
Assume that $T_\zeta\sim\exp(\zeta)$ and $\TAU_i \sim \exp(\beta+\omega)$, $i=1,2,\dots ~$. Then for $\beta , \omega > 0$ we have that
\begin{align}\label{eq:levy4}
\bar Y(T_\beta)&\sim \bar Y(T_{\beta+\omega})+W_{N_{\beta,\omega}} , \nonumber\\
W_{N_{\beta,\omega}}&\sim I\,(\bar Y(T_\beta)-V)^+ ,
\end{align}
where 
\begin{align}
    \label{defIW}
I&\sim 1_{\{N_{\beta,\omega}\ge 1\}}\sim {\mathbb B}(q), \nonumber\\
W_{N_{\beta,\omega}} &\sim \max_{k\in\{0,\ldots,N_{\beta,\omega}\}} Y(\Sigma_k),\nonumber\\
V&\sim \bar Y(T_{\beta+\omega})-Y(T_{\beta+\omega}),
\end{align}
and the random variables on the right of each of the equations in \eq{levy4} are assumed independent.
\end{Corollary}

{\it Proof.}
With $\TAU_i\sim \exp(\beta+\omega)$, it is easily verified that one has that $\Sigma_{N_{\beta,\omega}+1}\sim\exp(\beta)$.
In addition, $\Sigma_{N_{\beta,\omega}} \sim I\,\Sigma_{N_{\beta,\omega}+1}\sim I\,T_\beta$,
where $I,T_\beta$ are independent, so that
\begin{equation}
  (\bar Y(\Sigma_{N_{\beta,\omega}})-V)^+\sim I\,(\bar Y(T_\beta)-V)^+.
    \end{equation}
    The result now follows from
\cor{Levy2}.
\hfill$\Box$

\begin{Remark}\rm
Since for any L\'evy process $Y$ it holds that $\bar Y(t)-Y(t) \sim -\underline{Y}(t)$ where $\underline{Y}(t):=\inf_{s\in[0,t]}Y(s)$, we may replace what is written by $V_i\sim-\underline{Y}(\TAU_{i+1})$ or $V\sim -\underline{Y}(\TAU)$, in the obvious places. \hfill$\Diamond$
\end{Remark}
\noindent
Now we are in a position to claim that the above recursion approach answers Question 2 of Section~\ref{intro} affirmatively. Indeed, the simple recursions \eqref{eq:ZW} and \eqref{eq:ZW0} bring us via elementary steps to an alternative proof of Theorem~\ref{Thm1}, by making appropriate choices of $U_{n+1},V_n$ in \eqref{eq:ZW0}.
It should be observed that another elementary proof one might hope for is to rely on conditioning on the first event (being either an inspection or a killing), with the running maximum being $\bar Y(T_{\beta+\omega})$
until this event, but so far we have not been successful in proving Theorem~\ref{Thm1} along these lines.

Question 3 (concerning extensions of Theorem~\ref{Thm1}, while remaining in the one-dimensional setup) is also answered affirmatively.
Firstly, one has the transient decomposition in Corollary~\ref{cor:Levy}; secondly, other choices of $U,V$ than the one made in that corollary are possible; and thirdly, we mention three further generalizations in the remainder of this section.
\begin{Remark}\rm
A consequence of \cor{Levy3} is that we may repeat this decomposition for an increasing sequence of $\beta_n$ diverging to infinity (with $\beta_0=\beta$) and obtain (since $\bar Y(T_{\beta_n})$ vanishes in distribution) that
\begin{equation}\label{eq:R1}
\bar Y(T_\beta)\sim \sum_{n=1}^\infty R_n , 
\end{equation}
where $R_n$ are independent with \begin{equation}R_n\sim \max_{k\in\{0,\ldots, N(\beta_{n-1}/\beta_n)\}}S_k(\beta_n), \end{equation} and $S_k(\beta)$ is a random walk with i.i.d.\ increments distributed like $Y(T_\beta)$ and $N(p)+1\sim{\mathbb G}(p)$ is independent of the random walk. 
A special case is $\beta_n=\beta c^n$ where $c>1$, which gives
$R_n\sim \max_{n\in\{0,\ldots, N(1/c)\}}S_n(\beta c^n)$.

Letting $R_0\sim \sup_{k\in{\mathbb N}}S_k(\beta)$ we also have from \eq{levy2}:
\begin{equation}\label{eq:R2}
\bar Y(\infty) \sim\sum_{n=0}^\infty R_n\,.
\end{equation}
Note that the left hand sides of \eq{R1} and \eq{R2} (and hence also the right hand sides) are independent of the particular choice of the sequence $(\beta_n)_{n\in{\mathbb N}}$. For the special case $\beta_n=\beta c^n$, they are independent of the particular choice of $c$. 
\hfill$\Diamond$
\end{Remark}

\begin{Remark}\label{rem:RW}\rm
If instead of a L\'evy process, $(Y_n)_{n\in{\mathbb N}}$ is a random walk with i.i.d.\ increments, $\TAU,\TAU_1,\TAU_2,\ldots\sim{\mathbb G}(p)$ with $p\in(0,1)$ and $N+1\sim{\mathbb G}(p')$, then it is easy to check that $\Sigma_{N+1}\sim{\mathbb G}(pp')$; this is a consequence of the memoryless property (but one could alternatively prove this by using generating functions).
It is well known that, with $\bar Y_n=\max_{k\in\{0,\ldots,n\}}Y_k$, $\bar Y_{\TAU}$ and $\bar Y_{\TAU}-Y_{\TAU}$ (with $\TAU \sim {\mathbb G}(p)$) 
are independent in this case as well.
Hence, with  $T_\beta$ and $T_{\beta+\omega}$ being replaced by the random variables ${\mathbb G}(pp')$ and ${\mathbb G}(p)$, respectively, Corollaries~\ref{cor:Levy}-\ref{cor:Levy3} (and the remarks that follow) are valid for this discrete time case as well.

Reverting back to the L\'evy case, the above may be applied, in particular, to the random walk $(Y_{na})_{n\in{\mathbb N}}$ for any given $a>0$ (having infinitely divisible increments).   
\hfill$\Diamond$
\end{Remark}

Along the lines of the proof of \cor{Levy}, the following may also be easily shown.

\begin{Corollary}\label{cor:MAP}
Assume that $Y^1,Y^2,\ldots$ are L\'evy processes and that $\TAU_1,\TAU_2,\ldots$ are exponentially distributed. Neither sequence is assumed identically distributed. Assume that $Y^1,Y^2,\ldots,\TAU_1,\TAU_2,\ldots$ are independent. Denote $\Sigma_0:=0$ and $\Sigma_n:=\sum_{i=1}^n \TAU_i$ for $n\in\{ 1,2,\ldots\}$. For a given $n\in\{ 1,2,\ldots\}$ define the process $\{X^n(t),t\in[0,\Sigma_n]\}$ as follows:
\begin{equation}
X^n(t):=\begin{cases}
Y^n(t) , & 0\le t<\TAU_n , \\
Y^n(\TAU_n)+Y^{n-1}(t-\TAU_n),&\TAU_n\le t<\TAU_n+\TAU_{n-1},\\
\vdots&\vdots\\
Y^n(\TAU_n)+\ldots+Y^2(\TAU_2)+Y^1(t-(\Sigma_n-\TAU_1)),&\Sigma_n-\TAU_1\le t<\Sigma_n\,.
\end{cases}
\end{equation}
Finally, set $Z_n:=\sup_{s\in[0,\Sigma_n]}X^n(s)$, $U_i:=\sup_{s\in[0,\TAU_i]}Y^i(s)$ for $i\in\{1,2,\ldots\}$ and $V_i:=U_{i+1}-Y^{i+1}(\TAU_{i+1})$ for $i\in {\mathbb N}$.
Then
\begin{align}\label{eq:levy0}
Z_{n+1}\sim U_{n+1}+W_n, \nonumber\\
W_n\sim (Z_n-V_n)^+,
\end{align}
where, with $S_0:=0$ and $S_n:=\sum_{i=1}^n (U_i-V_i)$, $W_n=S_n-\min_{k\in\{0,\ldots, n\}}S_k$.
Moreover, the pairs of random variables appearing on the right of each equation in \eq{levy0} are independent.
\end{Corollary}

\section{Unique solution to a Lindley type equation}\label{Lindley}
In this section we study a generalization of the classical Lindley recursion, in which $W_{n+1} = (W_n+X_n)^+$, with special attention to its steady-state counterpart $W \sim (W+X)^+$.
In \prop{WXI} we provide a unique solution to the equation $W \sim I\,(W+X)^+$, with $I$ a Bernoulli random variable. We subsequently derive \cor{ZUVI}
from it, and observe that this corollary immediately implies \cor{Levy3} -- and hence Theorem~\ref{Thm1}.
At the end of the section we also discuss Question 4 posed in Section~\ref{intro}.

In preparation we first observe the following.
From the proof of \cor{Levy} it follows that if $T,X$ are nonnegative independent random variables then
\begin{equation}
\bar Y(T+X)\sim U+(\bar Y(T)-V)^+ ,
\end{equation}
where $\bar Y(T)$ and $(U,V)\sim (\bar Y(X),\bar Y(X)-Y(X))$ are independent. In this case $U,V$ are not necessarily independent.
\\
Now assume that $T$ is exponentially distributed. In addition, take $T'$ distributed like $T$ and independent of $T$ and $X$.
Using the memoryless property it is easily  checked that $T\sim X\wedge T'+1_{\{X<T'\}}T$ and we may replace $T,X$ in the previous paragraph by $1_{\{X<T'\}}T',X\wedge T'$ and obtain
\begin{equation}\label{eq:UVI}
\bar Y(T) \sim U+I\,(\bar Y(T)-V)^+ ,
\end{equation}
where \begin{equation}(U,V,I)\sim (\bar Y(X\wedge T),\,\bar Y(X\wedge T)-Y(X\wedge T),\,1_{\{X<T\}})\end{equation} is independent of $\bar Y(T)$.
One could also prove this by distinguishing between $X \geq T$, for which (\ref{eq:UVI}) trivially holds, and $X<T$, 
for which case one can follow the reasoning leading to (\ref{eq20}).
\begin{Remark}\label{rem:X}\rm
When $X$ is also exponentially distributed, then $X\wedge T$ and $I$ are independent, $X\wedge T$ is exponentially distributed and $U,V,I$ are independent.
Otherwise this is false. If ${\mathbb P}(X=0)<1$ (otherwise, $U=V=0$ and $I=1$), then $X\wedge T$ and $1_{\{X<T\}}$ are independent iff $X$ is also exponentially distributed.
To see this, note that if ${\mathbb P}(X=t)>0$ for some $t$ then in order for the events $\{\min(X,T)=t\}$ and $\{T>X\}$
to be independent it is necessary that ${\mathbb P}(T=t)>0$ which does not hold.
Therefore the distribution of $X$ must also be continuous and the assumptions of \cite{s65} are satisfied.
Using \cite{s65} we may now conclude that ${\mathbb P}(X>t)={\mathbb P}(T>t)^a$ for some $a>0$ and all $t$, where the right hand side, and hence also the left hand side, is the tail of an exponential distribution. $\hfill\Diamond$
\end{Remark}
In view of \eq{UVI} it is of interest to explore whether there is a unique distribution satisfying an equation of the form
\begin{equation}\label{eq:ZUVI}
Z\sim U+I\,(Z-V)^+,
\end{equation}
where $U,V\ge 0$, $I$ is an indicator and $U,I,V,Z$ on the right (the distributions of $U,V,I$ are considered known) are assumed independent. Consider $W\sim I(Z-V)^+$ when $Z$ satisfies \eq{ZUVI}. Then, taking $(U',V',I')$ to be an independent copy of $(U,V,I)$ which is independent of $Z$,
\begin{align}\label{eq:ZWUVI}
W&\sim I\,(Z-V)^+\sim I\,(U'+I'\,(Z-V')^+-V)^+\nonumber\\
&\sim I\,(U'+W-V)^+\sim I\,(W+U-V)^+\,,
\end{align}
where $U,V,I,W$ on the right are assumed independent. 
\begin{Proposition}\label{prop:WXI}
Assume that $I\sim{\mathbb B}(1-p)$, $p\in(0,1)$. Then
\begin{equation}\label{eq:IWX}
W\sim I\,(W+X)^+,
\end{equation}
where $X,I,W$ on the right are assumed independent, iff
\begin{equation}
W\sim \bar{S}_N\equiv\max_{k\in\{0,\ldots, N\}}S_k\,,
\end{equation} 
where $S\equiv (S_k)_{k\in{\mathbb N}}$ is a random walk with i.i.d.\ increments distributed like $X$, $N$ is independent of $S$ and $N+1\sim{\mathbb G}(p)$.
\end{Proposition}
We note that unlike Lindley's equation $W\sim (W+X)^+$, where one must assume that $S_n\to-\infty$ a.s., other than $p\in(0,1)$ no further stability requirements are needed here.

\medskip

{\it Proof.}
Assume that $W_0,X_1,X_2,\ldots,I_1,I_2,\ldots$ are independent. $P(W_0\ge 0)=1$ and for $i\ge 1$, $X_i\sim X$, $I_i\sim I$. Define the following recursion:
\begin{equation}
W_n=I_n(W_{n-1}+X_n)^+,
\end{equation}
for $n\in\{ 1,2,\ldots\}$. Then $(W_n)_{n\in{\mathbb N}}$ is a (possibly delayed) regenerative process with regeneration epoch $N+1=\inf\{n: I_n=0\}$, which has an aperiodic distribution.
In every regenerative cycle starting from the second, the process starts from zero. Therefore $W_n$ converges in distribution to some $W_\infty$ of which
the distribution is independent of the distribution of $W_0$. In order to identify this distribution, let us assume that $W_0=0$. In this case, $W_n=(W_{n-1}+X_n)^+$ for $n\in\{1,\ldots,N\}$, thus for each such $n$,
\begin{equation}
W_n=S_n-\min_{k\in\{0,\ldots, n\}}S_k\sim \max_{k\in\{0,\ldots, n\}}S_k\equiv \bar{S}_n\,,
\end{equation} 
where $S$ is a random walk with increments $X_1,X_2,\ldots$. Regenerative theory implies that, for every nonnegative (or bounded) Borel function $f$,
\begin{align}\label{eq:regen}
{\mathbb E}\,f(W_\infty)&=\frac{1}{{\mathbb E}(N+1)}{\mathbb E}\sum_{n=0}^N f(W_n)=p\,{\mathbb E}\sum_{n=0}^N f(W_n)\nonumber\\
&=p\,{\mathbb E}\sum_{n=0}^\infty f(W_n)1_{\{N\ge n\}}=\sum_{n=0}^\infty {\mathbb E}\,f(W_n)(1-p)^np\nonumber\\
&=\sum_{n=0}^\infty {\mathbb E}\,f(\bar{S}_n)(1-p)^np={\mathbb E}\,f(\bar{S}_N)\,,
\end{align}
so that $W_\infty\sim \bar{S}_N$.
Clearly $W_\infty\sim I\,(W_\infty+X)^+$ where $I,X,W_\infty$ on the right are independent, so that a solution to \eq{IWX} exists. For any distribution of $W$ that satisfies \eq{IWX}, we can start from $W_0\sim W$ and infer that $W_n\sim W$ for all $n\ge 0$ and thus also $W_\infty\sim W$. This implies uniqueness.
\hfill$\Box$

\medskip

In view of \eq{ZUVI} and \eq{ZWUVI}, the following is now immediate.
\begin{Corollary}\label{cor:ZUVI}
Assume that $U,V\ge 0$, and that $I\sim{\mathbb B}(1-p)$ for $p\in(0,1)$. Then
\begin{equation}
Z\sim U+I\,(Z-V)^+\,,
\end{equation}
where $U,V,I,Z$ on the right are assumed independent, iff
\begin{equation}
Z\sim U+W\,,
\end{equation}
where $W\sim\bar{S}_N$ is the solution from \prop{WXI} with $X\sim U-V$ and $U,W$ assumed independent. Moreover, $W\sim I\,(Z-V)^+$.
\end{Corollary}

We now observe that \cor{Levy3} ({\em and hence Theorem~\ref{Thm1}}) and similarly the related discrete time version (see \rem{RW}), can also be viewed as a consequence of \cor{ZUVI}. To see this,
recall that $T_\beta\sim\exp(\beta)$ for $\beta>0$. We have that $T_\beta\sim T_{\beta+\omega}+I\,T_\beta$ where $T_{\beta+\omega},T_\beta,I$ on the right are assumed independent and $I\sim{\mathbb B}(\omega/(\beta+\omega))$. Now,
\begin{equation}
\bar Y(T_\beta)\sim \bar Y(T_{\beta+\omega}+I\,T_\beta)\sim U+(\bar Y(I\,T_\beta)-V)^+=U+I\,(\bar Y(T_\beta)-V)^+\,,
\label{eq42}
\end{equation}
where $U\sim \bar Y(T_{\beta+\omega})$, $V\sim \bar Y(T_{\beta+\omega})-Y(T_{\beta+\omega})$, and $U,V,I,\bar Y(T_\beta)$ are independent. As $U-V\sim Y(T_{\beta+\omega})$
and thus \begin{equation}
    S_n\sim 
     Y(T_{\beta+\omega,1} + \dots + T_{\beta+\omega,n}),
\end{equation} \cor{ZUVI} implies \cor{Levy3}.

It is known that when $S_n\to-\infty$ a.s., there is a unique solution to $W\sim (W+X)^+$ where $W,X$ are independent which satisfies $W\sim\sup_{k\in{\mathbb N}}S_k$. Therefore in this case there is a unique solution to $Z\sim U+(Z-V)^+$ where $U,V\ge 0$ and where $U,V,Z$ on the right are independent, which satisfies $Z\sim U+W$ (as well as $W\sim (Z-V)^+$). This may be used to give an alternative derivation of \cor{Levy1}.

\begin{Remark}\label{rem:Ne}\rm
A similar approach may be applied to the recursion $W_n=I_n(W_{n-1}+X_n)^+$ where the $I_n$ are independent Bernoulli random variables, but not necessarily identically distributed. Whenever $I_n=0$ we probabilistically restart the sequence of indicators in order to obtain a natural regeneration epoch.
Namely, if $N$ is an integer valued random variable, then, whenever ${\mathbb P}(N\ge n-1)>0$,
let $I_n\sim{\mathbb B}(p_n)$, where $p_n={\mathbb P}(N\ge n)/{\mathbb P}(N\ge n-1)$. Then $N=\inf\{n: I_{n+1}=0\}$ and, similar as in \eq{regen}, we obtain that the limiting distribution of $W_n$ is the distribution of $\bar{S}_{N_e}$ where \begin{equation}{\mathbb P}(N_e=n):=\frac{{\mathbb P}(N\ge n)}{{\mathbb E}N+1}.\end{equation} \cor{ZUVI} may be modified accordingly. Note that in case $N+1\sim{\mathbb G}(p)$ we have that $N_e\sim N$. $\hfill\Diamond$
\end{Remark}

\begin{Remark}\rm 
\label{remarkques4}
We close this section discussing Question 4, raised in Section~\ref{intro}.
There we asked for an explanation of the intriguing fact that the inspection rate $\omega$ only appears in the right hand side of (\ref{onedimdeco}),
and not in its left hand side.
The same appears to hold in the two-dimensional decomposition of Theorem~\ref{Thm2}.
While it is intuitively obvious that the first vector in the right hand side of that decomposition is decreasing in $\omega$ and the second vector is increasing,
it is remarkable that the two $\omega$ contributions apparently cancel;
it indicates that there is an intimate relationship between $(\bar Y(T_{\beta+\omega}),G(T_{\beta+\omega}))$ and $(\bar S_{N_{\beta,\omega}},G_{N_{\beta,\omega}})$.
In the one-dimensional case of Theorem~\ref{Thm1}, the explanation of this phenomenon seems to lie in the distributional identity $T_\beta \sim T_{\beta+\omega} + I\,T_\beta$ with $I \sim {\mathbb B}(\omega/(\beta+\omega))$
and $T_{\beta+\omega}$, $T_\beta$ and $I$ independent. Observe that this identity has the same feature of $\omega$ only appearing in the right hand side.
This identity is used to derive (\ref{eq42}), which is an essential step in the alternative proof of Corollary~\ref{cor:Levy3} and hence of Theorem~\ref{Thm1}.
\hfill$\Diamond$
\end{Remark}

\section{A two-dimensional recursion}
\label{twodimrecur}
The goal of the present section is to give an elementary and transparent proof of the two-dimensional decomposition of Theorem~\ref{Thm2}.
We do so using a recursion of the type that was used in Section~\ref{onedimrecur} to deal with its one-dimensional counterpart (i.e., Theorem~\ref{Thm1}).

Theorem~\ref{Thm2} has, next to the running maximum, as a second component the time epoch at which that maximum last occurs.
Hence, while starting from a Lindley recursion between deterministic numbers $w_1,w_2,\dots,w_n$ (just like in Section~\ref{onedimrecur}), we now add as second component
the first index at which the largest of those $w_j$ occurs.
It should be noticed that there might be more than one index at which the largest $w_j$ occurs, which requires some care in using weak and strict inequalities; we will be working with the {\it first} index in order to, at a later stage, end up with relations for the {\it last} time epoch
(before $T_\beta$) at which the supremum of the L\'evy process $Y$ is attained. 
We pay much attention to the set-up of the recursion, but after the right recursion has been established,
we shall be more brief, as the procedure strongly resembles the one used in the one-dimensional case. Along the way we indicate in various places what the one-dimensional counterpart of a two-dimensional result is.

Let $x_1,x_2,\ldots$ a sequence of real numbers and $x'_n$ a sequence of nonnegative real numbers. Consider the following two-dimensional recursion with $(w_0,w'_0)=(0,0)$:
\begin{align}\label{eq:wprime}
w_n&:=(w_{n-1}+x_n)^+,\nonumber\\
w'_n&:=(w'_{n-1}+x'_n)1_{\{w_{n-1}+x_n \geq 0\}}\,.
\end{align}
Define the vectors $(s_0,s'_0):=(0,0)$ and $(s_n,s'_n):=\sum_{i=1}^n (x_i,x'_i)$ for $n\in\{ 1,2,\ldots\}$. With ${\mathbb I}_n:=\{0,1,\ldots,n\}$, let
\begin{equation}
j_n:=\min\left\{j\in {\mathbb I}_n:\,s_j=\min_{i\,\in\, {\mathbb I}_n}s_i\right\}\,.
\end{equation}
Then for every $j \in \{j_n,\ldots,n\}$, we have that $\min_{i\in\{1,\ldots, j\}}s_i=s_{j_n}$ and thus
\begin{equation}
w_j=s_j-\min_{i\in\{1,\ldots, j\}}s_i=s_j-s_{j_n}\,.
\end{equation}
Also, when $j_n<n$, then
we have that $s_j \geq s_{j_n}$ for $j\in\{j_n,\ldots, n\}$ and thus 
\begin{equation}w_{j-1}+x_j=s_{j-1}-s_{j_n}+x_j=s_j-s_{j_n} \geq 0.\end{equation}
When $j_n>0$ it clearly follows that $w_{j_n-1}+x_{j_n} <  w_{j_n}=0$ and thus, by \eq{wprime}, $w'_{j_n}=0$. When $j_n=0$ then $w'_0=0$ by assumption.
Observe that we had a strict inequality in the previous line, because $j_n$ is the {\em first} index at which the overall minimum is attained.
Recalling that $x'_1,x'_2,\ldots$ are nonnegative,  it follows that 
\begin{equation}
w'_n=s'_n-s'_{j_n}\,.
\end{equation}
Let us denote by $m'_n,m_n,k_n$ the values corresponding to $w'_n,w_n,n-j_n$ when the sequence $((x_1,x'_1),\ldots,(x_n,x'_n))$ is reversed, that is, replaced by $((x_n,x'_n),\ldots,(x_1,x'_1))$. Then, since we can also write
\begin{align}
w_n&=\max_{k\in {\mathbb I}_n}(s_n-s_{n-k}),\nonumber\\
n-j_n&=n-\min\left\{j\in {\mathbb I}_n|\,s_j=\min_{i\in {\mathbb I}_n}s_i\right\}=\max\left\{k\in {\mathbb I}_n|\,s_{n-k}=\min_{i\in {\mathbb I}_n}s_i\right\}\nonumber\\
&=\max\left\{k\in {\mathbb I}_n|\,s_{n-k}=\min_{i\in {\mathbb I}_n}s_{n-i}\right\}\nonumber\\ 
&=\max\left\{k\in {\mathbb I}_n|\,s_n-s_{n-k}=\max_{i\in {\mathbb I}_n}(s_n-s_{n-i})\right\},\nonumber\\
w'_n&=s'_n-s'_{n-(n-j_n)}\,,
\end{align}
it follows that
\begin{align}\label{eq:mk}
m_n&=\bar{s}_n=\max_{k\in {\mathbb I}_n}s_k=s_{k_n},\nonumber\\
k_n&=\max\{k\in {\mathbb I}_n|\, s_k=\bar{s}_n\},\nonumber\\
m'_n&=s'_{k_n}\,.
\end{align}
Hence, replacing $(x_i,x'_i)$ by i.i.d.\ random pairs $(X_i,X'_i)$ and similarly replacing $w_n,w_n',k_n,s_n$ by random variables with corresponding capital letters, we have that
\begin{equation}
(W_n,W'_n)\sim (\bar{S}_n, S'_{K_n})\,,
\end{equation}
where $K_n:=\max\{k\in {\mathbb I}_n:S_k=\bar{S}_n\}$ (the index of the last running maximum until $n$, that is) and we note that it is easy to verify that an equivalent definition of $K_n$ is $K_n:=\max\{k\in {\mathbb I}_n:S_k=\bar{S}_k\}$. The latter follows from the maximality of $K_n$, so that for $j\in\{K_n+1,\ldots,n\}$ (when $K_n<n$) we have that $S_j<\bar{S}_n$ while $\bar{S}_j=\bar{S}_n$ and thus $S_j<\bar{S}_j$.

Consider now the recursive equation with $(z_0,z'_0)=0$ and $u_n,u'_n,v_n,v'_n\ge 0$ for $n\in {\mathbb N}$ (cf.\ its one-dimensional counterpart in \eqref{eq:ZW}):
\begin{equation}
(z_{n+1},z'_{n+1})=(u_{n+1},u'_{n+1})+((z_n,z'_n)+(-v_n,v'_n))1_{\{z_n \geq v_n\}}\,.
\end{equation}
Then $(w_n,w'_n)=(z_{n+1}-u_{n+1},z'_{n+1}-u'_{n+1})$ satisfies \eq{wprime} with $x_n=u_n-v_n$ and $x'_n=u'_n+v'_n\ge 0$.
Therefore, when $((u_n,u'_n,v_n,v'_n))_{n\in{\mathbb N}}$ are replaced by i.i.d.\ random vectors, then this results in i.i.d.\ random pairs $(X_i,X'_i)$. As a consequence (cf.\ its one-dimensional counterpart in \prop{ind}),
\begin{equation}
(Z_{n+1},Z'_{n+1})\sim (U_{n+1},U'_{n+1})+(\bar{S}_n,S'_{K_n})\,.
\end{equation}  
Since the random pairs on the right are independent, we may take $(U,U')\sim (U_1,U'_1)$ which is independent of everything else, and conclude that
\begin{equation}\label{eq:Z2dim}
(Z_{n+1},Z'_{n+1})\sim (U,U')+(\bar{S}_n,S'_{K_n}),
\end{equation}  
where the two random pairs on the right are independent. As before, $n$ may be replaced by an independent nonnegative integer valued random variable $N$.
When $S_n\to-\infty$ a.s.\ ({\em e.g.,} when ${\mathbb E}X_1$ exists and is nonpositive), then $K_\infty,\bar{S}_\infty$ are a.s.\ finite and this also holds for random variables for which ${\mathbb P}(N=\infty)>0$.

We now point out how two-dimensional counterparts of Corollaries~\ref{cor:Levy1}, \ref{cor:Levy2} and \ref{cor:Levy3} -- and hence finally Theorem~\ref{Thm2} -- are obtained by an appropriate choice of $U,U',V,V'$.
For any $\tau,t\ge 0$,
\begin{align}
G(\tau+t)&=G(\tau)1_{\left\{Y(\tau)+\overline{\sigma_\tau Y}(t)<\bar{Y}(\tau)\right\}}+(\tau+\sigma_\tau G(t))1_{\left\{Y(\tau)+\overline{\sigma_\tau Y}(t)\ge\bar{Y}(\tau)\right\}}\nonumber\\
&=G(\tau)+ \Big(\sigma_{\tau}G(t)+(\tau-G(\tau)) \Big)
1_{\left\{\overline{\sigma_\tau Y}(t)\ge\bar{Y}(\tau)-Y(\tau)\right\}}\,,
\end{align}
while we observe that (\ref{eq22}) may be written as follows:
\begin{equation}
\bar{Y}(\tau+t)=\bar{Y}(\tau)+\left(\overline{\sigma_\tau Y}(t)-(\bar{Y}(\tau)-Y(\tau))\right)
1_{\left\{\overline{\sigma_\tau Y}(t)\ge\bar{Y}(\tau)-Y(\tau)\right\}}\,.
\end{equation}
Combining these two equations yields
\begin{align}
&(\bar{Y}(\tau+t),G(\tau+t))=(\bar{Y}(\tau),G(\tau))\nonumber\\
&\qquad+\left( (\overline{\sigma_\TAU Y}(t),\sigma_\tau G(t))+(-(\bar{Y}(\tau)-Y(\tau)),\tau-G(\tau))\right) 1_{\left\{\overline{\sigma_\tau Y}(t)\ge\bar{Y}(\tau)-Y(\tau)\right\}}\,,
\label{eq5353}
\end{align}
where $(\overline{\sigma_\tau Y}(t),\sigma_\tau G(t))$ is distributed like $(\bar{Y}(t),G(t))$ and is independent of everything else. 

Now replace $\tau,t$ by random variables. When $\tau$ is an exponentially distributed random variable then $(\bar{Y}(\tau),G(\tau))$ and $(\bar{Y}(\tau)-Y(\tau),\tau-G(\tau))$ are independent.
Taking $\tau, \tau_1,\tau_2,\ldots$ to be i.i.d. $\exp(\beta+ \omega)$ distributed and $t = \Sigma_n=\sum_{i=1}^n \tau_i$ then, as in Corollary~\ref{cor:Levy}
and via the derivation that led to \eq{Z2dim},
we end up with the equation
\begin{equation}
(\bar{Y}(\Sigma_{n+1}),G(\Sigma_{n+1})) \sim (\bar{Y}(T_{\beta+\omega}),G(T_{\beta+\omega}))+(\bar{S}_n,S'_{K_n}),
\end{equation} 
where
$S'_{K_n}$ is the time until we observe the last maximum of the random walk
with increments distributed like $Y(T_{\beta+\omega})$ until the $n$-th Poisson inspection.
When $n\to\infty$ we obtain the $(\bar{Y},G)$ version of \cor{Levy1} and for any distribution of $N$  we obtain this version of \cor{Levy2}.
In particular, when $N+1$ is geometrically distributed we have the $(\bar Y,G)$ version of \cor{Levy3} which implies Theorem~\ref{Thm2}.
\begin{Remark}\rm
We briefly return to Question 4 of Section~\ref{intro} regarding the fact that the inspection rate $\omega$ only appears on the right hand side of (\ref{onedimdeco})
and not the left hand side. We already observed in Remark~\ref{remarkques4} that
the same appears to hold in the two-dimensional decomposition of Theorem~\ref{Thm2}.
An explanation in the one-dimensional case was found in the distributional identity $T_\beta \sim T_{\beta+\omega} + I\,T_\beta$ with $I \sim {\mathbb B}(\omega/(\beta+\omega))$
and $T_{\beta+\omega}$, $T_\beta$ and $I$ independent. This identity has the same feature of $\omega$ only appearing in the right hand side, and
this identity is used to derive (\ref{eq42}), which is an essential step in the alternative proof of Corollary~\ref{cor:Levy3} and hence of Theorem~\ref{Thm1}.
In the present section, we could have given a slightly different proof of Theorem~\ref{Thm2} by using a 
two-dimensional (i.e., $(\bar{Y},G)$) version of (\ref{eq42}). That version is obtained by taking
$\tau=T_{\beta+\omega}$ and $t = I\,T_\beta$ in (\ref{eq5353}), so that $\tau+t = T_\beta$.
We thus obtain insight into the absence of the inspection rate $\omega$ in the left hand side of the two-dimensional decomposition.\hfill$\Diamond$
\end{Remark}

\begin{Remark}\rm
It is straightforward to observe that the approach in this section is also valid when $x'$ is a {\em vector} of any finite dimension.
Therefore it also directly follows that the following (modest) generalization is possible: if $J=\{J(t),t\ge 0\}$ is a multivariate subordinator which is independent of everything else, then
\begin{equation}
(\bar{Y}(T_\beta),J(G(T_\beta)))\sim (\bar{Y}(T_{\beta+\omega}),J(G(T_{\beta+\omega})))+(\bar{S}_N,J(S'_{K_N}))\,,
\end{equation}
where the two pairs on the right are independent.\hfill$\Diamond$
\end{Remark}

Finally, \prop{WXI} (regarding the unique solution of the generalized Lindley equation) may be easily extended to the current case as well with the same proof. Namely, the unique solution of the equation
\begin{equation}
(Z,Z')\sim (U,U')+I\,((Z,Z')+(-V,V'))\,,
\end{equation}
where $I$ and the three pairs on the right are independent, is
\begin{equation}
(Z,Z')\sim (U,U')+(\bar{S}_N,S'_{K_N}),
\label{eqzz}
\end{equation}
where $N+1$ is geometrically distributed. 

Also, as in \rem{Ne}, considering instead a regenerative sequence with $N$  being an integer valued random variable with finite mean, we replace $N$ with $N_e$ in the right hand side of (\ref{eqzz}).

\section{Spectrally one-sided cases}
\label{implic}
In this section we briefly discuss the case that the L\'evy process $Y$ is spectrally positive or spectrally negative. As it turns out, we can then easily obtain explicit expressions for the joint transforms of $(\bar Y(T_\beta), G(T_\beta))$ and, hence, of $(\bar S_{N_{\beta,\omega}},G_{N_{\beta,\omega}})$.
The spectrally positive case is considered in Subsection~\ref{sec6.1}, and the spectrally negative case in Subsection~\ref{sec6.2}.
In both cases, the starting point is the expression given in Remark \ref{R0}. 

\subsection{The spectrally positive case}
\label{sec6.1}
When the L\'evy process $Y$ is spectrally positive, i.e., it has no downward jumps, then it follows from  \cite[Section 6.5.2]{kyprianou} that
\begin{equation}
\kappa(\beta,\alpha) = \frac{\beta-\varphi(\alpha)}{\psi(\beta)-\alpha},
\end{equation}
where $\varphi(\alpha)$ denotes the Laplace exponent $\log {\mathbb E} \exp(-\alpha Y(1))$ of the L\'evy process $Y$ and $\psi(\alpha)$ its right-inverse.
Hence it follows from \eqref{doubleLT} that
\begin{equation}
{\mathbb E} {\rm e}^{-\alpha \bar S_{N_{\beta,\omega}} - \gamma G_{N_{\beta,\omega}}} = \frac{\beta}{\psi(\beta)} \frac{\psi(\beta+\gamma)-\alpha}{\beta+\gamma-\varphi(\alpha)} \frac{\psi(\beta+\omega)}{\beta+\omega} \frac{\beta+\omega+\gamma-\varphi(\alpha)}{\psi(\beta+\omega+\gamma)-\alpha}.
\label{doubleLTSP}
\end{equation}
The known bivariate transform (cf.\ (\ref{WH1})) \begin{equation}\frac{\beta}{\psi(\beta)} \frac{\psi(\beta+\gamma)-\alpha}{\beta+\gamma-\varphi(\alpha)}
\end{equation}
for the continuously inspected case is recovered when sending $\omega$ to $\infty$. 
Furthermore, moments are readily obtained; in particular,
\begin{align}
    {\mathbb E} \bar{S}_{N_{\beta,\omega}} &= \frac{1}{\psi(\beta)} - \frac{1}{\psi(\beta+\omega)} - \frac{\varphi'(0)}{\beta} + \frac{\varphi'(0)}{\beta+\omega},\\
    {\mathbb E} G_{N_{\beta,\omega}} &= -\frac{\psi'(\beta)}{\psi(\beta)} + \frac{\psi'(\beta+\omega)}{\psi(\beta+\omega)} + \frac{1}{\beta} - \frac{1}{\beta+\omega},\\
    {\mathbb E}\left[\bar{S}_{N_{\beta,\omega}} G_{N_{\beta,\omega}}\right] &= 
   \frac{\psi'(\beta)}{\psi(\beta) \psi(\beta+\omega)} + \frac{\psi'(\beta+\omega)}{\psi(\beta)\psi(\beta+\omega)} - 2 \frac{\psi'(\beta+\omega)}{(\psi(\beta+\omega))^2}\:+
   \nonumber\\&\:\:\:\:\:\:
    \varphi'(0) \left(\frac{1}{\beta} - \frac{1}{\beta+\omega}\right)\left(\frac{\psi'(\beta)}{\psi(\beta)} - \frac{\psi'(\beta+\omega)}{\psi(\beta+\omega)}\right)\:+
   \nonumber\\
   &\:\:\:\:\:\: \left(\frac{1}{\beta} - \frac{1}{\beta+\omega}\right)\left(\frac{1}{\psi(\beta+\omega)}-\frac{1}{\psi(\beta)}\right) - 2\frac{\varphi'(0)}{\beta} (\frac{1}{\beta} - \frac{1}{\beta+\omega}) .
\end{align}
Hence
\begin{equation}
    {\mathbb C}{\rm ov}(\bar{S}_{N_{\beta,\omega}}, G_{N_{\beta,\omega}}) = - \varphi'(0)\left(\frac{1}{\beta^2} - \frac{1}{(\beta+\omega)^2} \right) + \frac{\psi'(\beta)}{(\psi(\beta))^2} - \frac{\psi'(\beta+\omega)}{(\psi(\beta+\omega))^2} . 
\end{equation}
Furthermore, assuming that $\varphi''(0) < \infty$,
\begin{equation}
    {\rm Var} \,\bar{S}_{N_{\beta,\omega}} = \varphi''(0) \left(\frac{1}{\beta} - \frac{1}{\beta+\omega}\right) + \varphi'(0))^2\left(\frac{1}{\beta^2} - \frac{1}{(\beta+\omega)^2}\right) - \left(\frac{1}{\psi(\beta)}\right)^2 + \left(\frac{1}{\psi(\beta+\omega)}\right)^2 , 
    \label{VarSbar}
\end{equation}
\begin{equation}
    {\rm Var} \,G_{N_{\beta,\omega}} = \frac{1}{\beta^2} - \frac{1}{(\beta+\omega)^2} + \frac{\psi''(\beta)}{\psi(\beta)} - \frac{\psi''(\beta+\omega)}{\psi(\beta+\omega)} - \left(\frac{\psi'(\beta)}{\psi(\beta)}\right)^2 + \left(\frac{\psi'(\beta+\omega)}{\psi(\beta+\omega)}\right)^2 .
    \label{VarG}
\end{equation}
Of course the correlation coefficient between the running maximum at inspection epochs before $T_\beta$ and the last epoch of its occurrence immediately follows from the last three formulas. Note that the terms in (\ref{VarSbar}) and (\ref{VarG}) involving $\beta$ (resp.\ $\beta+\omega$) reveal the variance of $\bar{Y}(T_\beta)$ and $\bar{G}(T_\beta)$ (resp.\ the variance of $\bar{Y}(T_{\beta+\omega})$ and $\bar{G}(T_{\beta+\omega})$).

\subsection{The spectrally negative case}
\label{sec6.2}
When the L\'evy process $Y$ is spectrally negative, i.e., it has no upward jumps, then it follows from \cite[Section 6.5.2]{kyprianou} that
\begin{equation}
\kappa(\beta,\alpha) = \Psi(\beta)+\alpha,
\end{equation}
where $\Psi(\alpha)$ is the right-inverse of the cumulant $\Phi(\alpha)= \log {\mathbb E} \exp(\alpha Y(1))$.
Hence it follows from \eqref{doubleLT} that
\begin{equation}
{\mathbb E} {\rm e}^{-\alpha \bar S_{N_{\beta,\omega}} - \gamma G_{N_{\beta,\omega}}} = \frac{\Psi(\beta)}{\Psi(\beta+\omega)} \frac{\Psi(\beta+\omega+\gamma)+\alpha}{\Psi(\beta+\gamma)+\alpha}.
\label{doubleLTSN}
\end{equation}
This reveals that $\bar S_{N_{\beta,\omega}}$ has an atom $\Psi(\beta)/\Psi(\beta+\omega)$ at zero and is exp($\Psi(\beta)$) distributed with the complementary probability. Also, when considering the limit as $\omega\to\infty$ we recover the known expression $\Psi(\beta)/(\Psi(\beta+\gamma)+\alpha)$.
Furthermore, the transform of the time of the last running maximum has the following elegant expression:
\[
{\mathbb E} {\rm e}^{- \gamma G_{N_{\beta,\omega}}} = \frac{\Psi(\beta)}{\Psi(\beta+\omega)} \frac{\Psi(\beta+\omega+\gamma)}{\Psi(\beta+\gamma)}.
\]
It requires some elementary calculus to verify that
\begin{align}
    {\mathbb E} \bar{S}_{N_{\beta,\omega}} &= \frac{1}{\Psi(\beta)} - \frac{1}{\Psi(\beta+\omega)},
   \\
{\mathbb E} G_{N_{\beta,\omega}} &= \frac{\Psi'(\beta)}{\Psi(\beta)} - \frac{\Psi'(\beta+\omega)}{\Psi(\beta+\omega)},  \\
    {\mathbb E} \left[\bar{S}_{N_{\beta,\omega}}G_{N_{\beta,\omega}}\right] &= 
    \frac{\Psi'(\beta)}{\Psi(\beta)} \left( \frac{1}{\Psi(\beta)} - \frac{1}{\Psi(\beta+\omega)}\right)
    + \frac{1}{\Psi(\beta)} \left( \frac{\Psi'(\beta)}{\Psi(\beta)} - \frac{\Psi'(\beta+\omega)}{\Psi(\beta+\omega)}\right),
\end{align}
so that
\begin{equation}
    {\mathbb C}{\rm ov}(\bar{S}_{N_{\beta,\omega}} ,G_{N_{\beta,\omega}}) = \frac{\Psi'(\beta)}{(\Psi(\beta))^2} - \frac{\Psi'(\beta+\omega)}{(\Psi(\beta+\omega))^2}.
\end{equation}
Furthermore,
\begin{align}
    {\mathbb V}{\rm ar} \,\bar{S}_{N_{\beta,\omega}}& = \frac{1}{(\Psi(\beta))^2} - \frac{1}{(\Psi(\beta+\omega))^2},\\
    {\mathbb V}{\rm ar} \,G_{N_{\beta,\omega}} &= \frac{\Psi''(\beta+\omega)}{\Psi(\beta+\omega)} - \frac{\Psi''(\beta)}{\Psi(\beta)} + \left(\frac{\Psi'(\beta)}{\Psi(\beta)}\right)^2 - \left(\frac{\Psi'(\beta+\omega)}{\Psi(\beta+\omega)}\right)^2 .
\end{align}

\section{Concluding remarks}
\label{sec:concl}
In this concluding section we discuss a number of ramifications of our work as well as directions for follow-up research.
\begin{itemize}
\item[$\circ$]
Our work has an interesting implication when studying, in the context of in insurance risk, {\it Parisian ruin}. Parisian ruin occurs for an insurance company when its capital level 
has been negative for a period of at least length $F$.
Albrecher and Ivanovs \cite{ai17} point out that, if lengths of successive excursions below zero are compared with i.i.d.\ random variables $F_1,F_2,\ldots \sim {\rm exp}(\omega)$ \cite{landriault}, then the ruin time for Parisian ruin
has the same distribution as the time until bankruptcy in the case of Poisson($\omega$) inspections.
Note that this observation adds to the importance of studying the time until bankruptcy.
Landriault et al.\ \cite{landriault} derive explicit expressions for the Laplace transform, with $P_x$ denoting the Parisian ruin time given that the initial capital was $x$, \begin{equation}
    g_\beta(x) := \int_0^{\infty} {\rm exp}(-\beta y) ~{\mathbb P}(P_x \in{\rm d}y)
\end{equation} in \cite[Corollary 3.1]{landriault}, and for the special case $x=0$ in \cite[Corollary 3.2]{landriault}. Notice that $g_\beta(x)$ can be viewed as the probability that Parisian ruin occurs before $T_\beta$. Taking the Laplace-Stieltjes transform of $-g_\beta(x)$ with respect to $x$ (realizing that $g_\beta(x)$ can be interpreted as the bankruptcy probability ${\mathbb P}(\bar{S}_{N_{\beta,\omega}} > x)$) and adding the contribution $g_\beta(0)$ due to the jump in zero, one eventually obtains the expression for ${\mathbb E} \exp(-\alpha \bar{S}_{N_{\beta,\omega}})$ given in (\ref{doubleLTSP}) (choosing $\gamma = 0$), as it should be. 

\item[$\circ$]
The {\it asymptotics} of the all-time bankruptcy probability in the spectrally one-sided case were identified in \cite{bm21}. This left open the case of the driving L\'evy process being spectrally two-sided. Recognizing that the process at inspection epochs is a random walk, in principle the results in \cite{korshunov} reveal these asymptotics, albeit in a rather implicit form. In \cite{korshunov} three cases were distinguished: the increments of the random walk having light tails, heavy tails, and being in an intermediate regime. It would be interesting to compare these results with their continuous-inspection counterparts \cite{bertoindoney,dieker}, so as to assess the `information loss' due to the Poisson inspection scheme. In addition, the joint asymptotics of $\bar{S}_{N_{\beta,\omega}}$ and $G_{N_{\beta,\omega}}$ could be an interesting research topic. 
\item[$\circ$]
In the present paper we have worked with Poisson inspection epochs, but one may wonder what happens if we work with {\it alternative inspection schemes}. It is conceivable that e.g.\ the case of Erlang inter-inspection times is amenable for analysis (cf.\ \cite[Remark 3.1]{bm21}), but one would ideally allow the more general class of renewal inspection processes. The results of \cite{korshunov} reveal that, {for renewal inter-inspection times and the driving L\'evy process $Y$ being light-tailed}, the asymptotics of the all-time bankruptcy probability are of the type $\gamma {\rm e}^{-\theta^\star u}$, with only the $\gamma$ depending on the distribution of the inter-inspection times (i.e., the decay rate $\theta^\star$ is unaffected). {In the case that the L\'evy process is heavy-tailed}, the results of \cite[Section 4.2]{bm21} suggest that the inspection process may not have any impact on the asymptotics of the bankruptcy probability. 
\item[$\circ$]

Importantly, our decomposition results do not provide us with insight into the {\it joint distribution} of $\bar Y(T_\beta)$ and $\bar S_{N_{\beta,\omega}}$, nor into the conditional distribution of $\bar Y(T_\beta)$ given the value of $\bar S_{N_{\beta,\omega}}$. For practical purposes, it would be useful to understand these relations, because that helps shedding light on the likelihood of having exceeded a certain threshold based on partial information. Ideally, one would be able to somehow get a handle on probabilities of the type
\begin{equation}{\mathbb P}
(\bar Y(T_\beta)>u \,|\,S_0,S_1,\ldots, {S}_{N_{\beta,\omega}})
\end{equation}
or 
\begin{equation}{\mathbb P}
(\bar Y(T_\beta)>u \,|\,\bar S_{N_{\beta,\omega}}).
\end{equation}
\end{itemize}

\end{document}